\documentclass[12pt]{article}
\usepackage{graphics}
\parskip=\medskipamount
\begin{document}
\title{Reliable Operations on Oscillatory Functions}
\author{{Gh.~Adam\thanks{~E-mail: adamg@theor1.theory.nipne.ro}\ 
~and
S.~Adam\thanks{~E-mail: adams@theor1.theory.nipne.ro}}
\\
{\em Department of Theoretical Physics,}\\
{\em Institute of Physics and Nuclear Engineering,}\\
{\em P.O.~Box MG-6, 76900 Bucharest-M\u agurele, Romania}\\
}
\maketitle

\begin{abstract}
Approximate $p$-point Leibniz derivation formulas as well as interpolatory
Simpson quadrature sums adapted to oscillatory functions are discussed.
Both theoretical considerations and numerical evidence concerning
the dependence of the discretization errors on the frequency parameter of
the oscillatory functions show that the accuracy gain of the present formulas
over those based on the exponential fitting approach [L.~Ixaru, {\sl Computer
Physics Communications}, 105 (1997) 1--19] is overwhelming.
\end{abstract}

\section{Introduction}
The mathematical descriptions of classical oscillatory phenomena (like
vibrations, wave propagation, resonances) or of the behaviour of quantum
systems involve operations on oscillatory functions. In most instances,
such operations (e.g., differentiation, integration, solving differential
equations) require numerical methods.

A successful approach towards accurate approximation of the oscillatory
solutions of ordinary differential equations is the exponential fitting
method, first proposed for the radial Schr\" odinger equation
%                     \cite Raptis [1978]
\cite{Raptis78},
and then extended by several authors to more general differential equations
(see, e.g.,~%
%                     \cite Ixaru [1997]
\cite{Ixaru97}
for a list of relevant results).

In a recent paper, Ixaru
%                     \cite Ixaru [1997]
\cite{Ixaru97}
raised the question whether the exponential fitting approach could yield
useful formulas for the numerical derivation and
integration of oscillatory functions as well. His investigation resulted in
several new formulas which extended to oscillatory functions well known
elementary discrete approximations of the first and second order
derivatives as well as the Simpson quadrature formula.

With the aim at enhancing code robustness and reliability over integration
subranges characterized by the occurrence of a slowly varying regular factor,
we tried to incorporate the optimized Simpson quadrature formula of 
%                     \cite Ixaru [1997]
\cite{Ixaru97}
in our computer code devoted to the integration of the products of functions
with oscillatory factors
%                     \cite Adam [1991]
%                     \cite Adam [1993]
%                     \cite Adam [1993]
\cite{Adam91,Adam93,Adam99}
within a class conscious automatic adaptive quadrature frame
%                     \cite DavRab [1984]
\cite{DavRab84}.
Since the attempt resulted in significant deterioration of code performances,
we undertook a complementary study the results of which are reported below.

There are four fundamental functional dependences of interest, namely,
\begin{equation}
   F(x) \equiv F_{s,\eta}^{\omega,\delta}(x)
          = f(x) g_{s,\eta}(\omega x + \delta)\,,
\label{ygeneric}
\end{equation}
where $f(x)$ is a sufficiently smooth real function, while
$g_{s,\eta}$ ($s = 1,2$ , $\eta = \pm 1$) denotes one of the following
four weight functions
\begin{equation}
\begin{array}{rclcrcl}
g_{1,-1}(\omega x+\delta) &\!\! = \!\!& \cos (\omega x+\delta)&\!\! , &
g_{2,-1}(\omega x+\delta) &\!\! = \!\!& \sin (\omega x+\delta)\,,\\
g_{1,1}(\omega x+\delta)  &\!\! = \!\!& \cosh (\omega x+\delta)&\!\! , &
g_{2,1}(\omega x+\delta)  &\!\! = \!\!& \sinh (\omega x+\delta)\,.
\end{array}
\label{ggeneric}
\end{equation}
Here, the frequency parameter $\omega$ and the initial phase $\delta$ are
both real and constant.

The principle of the present approach towards numerical differentiation and
numerical integration of functions of the form~(\ref{ygeneric}) or of their
linear combinations is well-known (see, e.g., the derivation of extended
Clenshaw-Curtis quadrature sums for oscillatory functions in~%
%                     \cite QUADPACK [1983]
\cite{QUADPACK}):
the approximating operations apply to the regular factor $f(x)$ only, while
the contribution of the weight factor $g_{s,\eta}$ is included {\sl exactly}. 
This results in formulas which are {\sl uniformly valid\/} at values of
the frequency parameter $\omega$ running over the whole set of machine
numbers at which the pair of weight functions \{$g_{1,\eta}(\omega x+\delta)$,
$g_{2,\eta}(\omega x+\delta)$\} can be accurately computed.
Of course, the classical differentiation and integration formulas are
recovered in the limit $\omega \rightarrow 0$.

The derivation of exponential fitting formulas of numerical differentiation
and integration of~%
%                     \cite Ixaru [1997]
\cite{Ixaru97},
requires the {\sl preservation of the formal structure of the corresponding
classical formulas}. The consequence is that the exponentially fitted
coefficients depend on the $g_{s,\eta}$ factor in a manner which results
in {\sl breakdowns of the obtained formulas\/} at some specific frequency
parameter values.
Due to this feature, the formulas of numerical differentiation and integration
based on exponential fitting are of limited practical value. In particular,
this is the reason for the abovementioned performance deterioration of the
code of automatic quadrature under incorporation of the optimal exponential
fitting generalization of Simpson quadrature formula.

To allow straightforward comparison of the $\omega$ dependences of the
discretization errors associated to the approximating formulas within the two
approaches, the case study problems considered in
%                     \cite Ixaru [1997]
\cite{Ixaru97}
have been solved by both methods. The accuracy improvement brought by the
present formulas with respect to those based on exponential fitting was found
to be comparable to that reported in~%
%                     \cite Ixaru [1997]
\cite{Ixaru97}
to be brought by the exponential fitting ones with respect to their classical
counterparts.

The paper is organized in three sections. Section~\ref{sec:deriv} discusses
discretized Leibniz derivation formulas. Section~\ref{sec:quad} provides
three-point interpolatory quadrature sums for integrands of the
form~(\ref{ygeneric}). Section~\ref{sec:concl} summarizes the main practical
consequences of the present investigation.
\section{Derivatives}
\label{sec:deriv}
\subsection{Discretized Leibniz formulas}
The description of oscillatory phenomena generally involves series of
harmonics, with terms of the form
\begin{equation}
  \Phi (x) = f_1(x) g_{1,\eta}(\omega x + \delta) +
         f_2(x) g_{2,\eta}(\omega x + \delta)\,.
\label{reffunc}
\end{equation}
To get numerical differentiation or integration formulas of such an
expression, solutions for the basic functional dependences~(\ref{ygeneric})
are needed.

Straightforward use of the Leibniz derivation formula yields for the
$n$-th derivative of (\ref{ygeneric})
\begin{equation}
   F^{(n)}(x) = \sum_{k=0}^{n}
                   \left(
                       \begin{array}{cc}
                         n\\
                         k
                        \end{array}
                   \right)
                     \omega^{k}f^{(n-k)}(x) g_{s,\eta}^{(k)}(t)\,,
\label{Leibniz}
\end{equation}
where, for each involved function $y(u)$,
$y^{(m)}(u) \equiv d^{m}y(u)/du^{m}$, while $t = \omega x + \delta \,.$
For any natural number $k$,
\begin{equation}
 \begin{array}{rclcrccc}
  g_{s,\eta}^{(4k+1)}(t) &\!\! = \!\!& \eta ^{s}  &\!\! g_{3-s,\eta}(t) \,, &
  g_{s,\eta}^{(4k+2)}(t) &\!\! = \!\!& \eta       &\!\!  g_{  s,\eta}(t)\,,\\
  g_{s,\eta}^{(4k+3)}(t) &\!\! = \!\!& \eta ^{s-1}&\!\! g_{3-s,\eta}(t) \,, &
  g_{s,\eta}^{(4k+4)}(t) &\!\! = \!\!&            &\!\!  g_{  s,\eta}(t)\,,
 \end{array}
\end{equation}
hence (\ref{Leibniz}) consists of a superposition of the pair of functions
$g_{1,\eta}(t)$ and $g_{2,\eta}(t)$, with coefficients which are themselves
expressed as superpositions of derivatives of the regular factor $f(x)$.

The use of the Leibniz formula (\ref{Leibniz}), combined with
classical $p$-point difference classical discretization formulas
for the involved derivatives of the regular factors $f_1(x)$ and $f_2(x)$
respectively will result in approximate derivation formulas of $\Phi (x)$,
Eq.~(\ref{reffunc}), which are {\sl optimal\/} with respect to the frequency
parameter~$\omega$. The $p$-point discretizations of the first and second
order derivatives of~$\Phi (x)$ are considered in detail below.
\subsubsection{First Order Derivative}
%\bigskip
%
%\noindent
%{\sl First Order Derivative}

The above equations yield
\begin{eqnarray}
  \Phi^{\prime}(x) &\!\! = \!\!& [f^{\prime}_1(x) + \omega f_2(x)]
                    g_{1,\eta}(\omega x + \delta)\nonumber \\
                &\!\! + \!\!& [f^{\prime}_2(x) + {\eta} \omega f_1(x)]
                    g_{2,\eta}(\omega x + \delta)\,.
\label{ypgeneric}
\end{eqnarray}

The present two-point approximation to $\Phi^{\prime}(x)$ which corresponds
to the exponentially fitted formula~(2.28) of
%                     \cite Ixaru [1997]
\cite{Ixaru97}
reads
\begin{eqnarray}
  \Phi^{\prime}_{2}(x) &\!\! = \!\!& \frac{1}{2h}\{[f_1(x+h) - f_1(x-h) +
                    2\lambda f_2(x)]
                    g_{1,\eta}(\omega x + \delta)\nonumber \\
                &   & \ +~[f_2(x+h) - f_2(x-h) + 2{\eta} \lambda f_1(x)]
                    g_{2,\eta}(\omega x + \delta)\}\,,
\label{ypoh2}
\end{eqnarray}
where
\begin{equation}
  \lambda = h \omega\;.
\label{deflambda}
\end{equation}
If the non-oscillatory factors $f_1(x)$ and $f_2(x)$ are at least three times
differentiable, then the discretization error associated to~(\ref{ypoh2})
{\sl shows ${\cal O} (h^2)$ accuracy}
\begin{eqnarray}
  e^{\prime}_2(x) &\!\! = &\!\! - \frac{1}{6}h^2[f_1^{(3)}(x+\theta _1 h)
                    g_{1,\eta}(\omega x + \delta)
        + f_2^{(3)}(x+\theta _2 h) g_{2,\eta}(\omega x + \delta)]\,,\nonumber\\
               & & {\;\; \; \; \; \; \; \; \; \; \; \; \; \; \; \; \; \;}
                -1 < \theta _1, \theta _2 < 1 \,.
\label{erh2}
\end{eqnarray}

Similarly, the present four-point approximation to $\Phi^{\prime}(x)$ which
corresponds to the exponentially fitted formula~(3.1) with coefficients~(3.3) of
%                     \cite Ixaru [1997]
\cite{Ixaru97}
reads
\begin{eqnarray}
  \Phi^{\prime}_{4}(x) &\!\! = \!\!& \frac{1}{12h}
            \{[f_1(x-2h) - 8 f_1(x-h) + 8 f_1(x+h)\nonumber\\
                & & {\; \; \; \; \; \; \; \; \;} - f_1(x+2h) +
                    12 \lambda f_2(x)] g_{1,\eta}(\omega x + \delta) \nonumber\\
                & & \ +~[f_2(x-2h) - 8 f_2(x-h) + 8 f_2(x+h)\nonumber\\
                & & {\; \; \; \; \; \; \; \; \;} - f_2(x+2h) +
                    12 {\eta} \lambda f_1(x)] g_{2,\eta}(\omega x + \delta)\}\,.
\label{ypoh4}
\end{eqnarray}
This expression holds for non-oscillatory factors $f_1(x)$ and $f_2(x)$
which are at least five times differentiable and associate an
{\sl ${\cal O} (h^4)$ discretization error}
\begin{eqnarray}
  e^{\prime}_4(x) &\!\! = &\!\! \frac{1}{30}h^4[f_1^{(5)}(x+\theta _3 h)
                    g_{1,\eta}(\omega x + \delta)
        + f_2^{(5)}(x+\theta _4 h) g_{2,\eta}(\omega x + \delta)]\,,\nonumber\\
               & & {\;\; \; \; \; \; \; \; \; \; \; \; \; \; \; \; \; \;}
                -1 < \theta _3, \theta _4 < 1 \,.
\label{erh4}
\end{eqnarray}
\subsubsection{Second Order Derivative}
%\bigskip
%
%\noindent
%{\sl Second Order Derivative}
%
For the second order derivative of (\ref{reffunc}), we get
\begin{eqnarray}
  \Phi^{\prime\prime}(x) &\!\! = \!\!&
      [f^{\prime\prime}_1(x) + 2\omega f^{\prime}_2(x)
       +\eta \omega ^2 f_1(x)]
                    g_{1,\eta}(\omega x + \delta)\nonumber \\
                      & + &
      [f^{\prime\prime}_2(x) + 2\eta \omega f^{\prime}_1(x)
       +\eta \omega ^2 f_2(x)]
                    g_{2,\eta}(\omega x + \delta)\,.
\label{ysgeneric}
\end{eqnarray}
hence the approximation to $\Phi^{\prime\prime}(x)$ which corresponds to the
the exponentially fitted formula (3.16) with coefficients (3.18) of
%                     \cite Ixaru [1997]
\cite{Ixaru97}
is provided by the expression
\begin{eqnarray}
  \Phi^{\prime\prime}_{2}(x) &\!\! = \!\!& 
    \frac{1}{h^2}
     \Bigl\{
       \{[f_1(x+h) + (\eta \lambda^2 - 2)f_1(x) + f_1(x-h)]
             \nonumber\\
          & & {\;} + \lambda [f_2(x+h) - f_2(x-h)]\}
               g_{1,\eta}(\omega x + \delta) \nonumber\\
   &\!\!+\!\!& \{[f_2(x+h) + (\eta \lambda^2 - 2)f_2(x) + f_2(x-h)] \nonumber\\
          & &
            {\;} + \eta \lambda [f_1(x+h) - f_1(x-h)]\}
               g_{2,\eta}(\omega x + \delta)
     \Bigr\}\,.
\label{ysoh2}
\end{eqnarray}
This expression is useful provided the non-oscillatory factors $f_1(x)$ and
$f_2(x)$ are at least four times differentiable. Then its associated
discretization error {\sl shows ${\cal O} (h^2)$ accuracy}
\begin{eqnarray}
  e^{\prime\prime}_2(x) &\!\! = \!\!& - \frac{1}{12}h^2
      \{[f_1^{(4)}(x+\theta _5 h) + 4 \omega f_2^{(3)}(x)]
              g_{1,\eta}(\omega x + \delta)\nonumber\\
  & & + [f_2^{(4)}(x+\theta _6 h) + 4 \eta \omega f_1^{(3)}(x)]
              g_{2,\eta}(\omega x + \delta)
      \}\,, \; \; -1 < \theta _5, \theta _6 < 1\,.
\label{er2h2}
\end{eqnarray}

{\sl The dependence of the discretization errors on the frequency parameter
$\omega$} will be discussed in the case of trigonometric functions
($\eta=-1$) only. Then $\vert \cos (\omega x + \delta )\vert \leq 1$ and
$\vert \sin (\omega x + \delta )\vert \leq 1$, hence Eqs.~(\ref{erh4})
and~(\ref{ysoh2}) result in the following {\sl upper bounds}, which are
{\sl independent of $\omega$ in the case of first order derivatives\/} and
{\sl linearly increasing with $\vert \omega \vert $ in the case of second
order derivatives\/}:
\begin{eqnarray}
  M^{\prime}_2(x) &\!\! = &\!\! 
      \frac{1}{6}h^2[\vert f_1^{(3)}(x+\theta _1 h)\vert 
                   + \vert f_2^{(3)}(x+\theta _2 h)\vert ]\,,
\label{mph2}\\
  M^{\prime}_4(x) &\!\! = &\!\! 
      \frac{1}{30}h^4[\vert f_1^{(5)}(x+\theta _3 h)\vert 
        + \vert f_2^{(5)}(x+\theta _4 h)\vert ]\,,
\label{mph4}\\
  M^{\prime\prime}_2(x) &\!\! = \!\!& \frac{1}{12}h^2
      \{[\vert f_1^{(4)}(x+\theta _5 h)\vert +
      4 \vert \omega \vert \vert f_2^{(3)}(x)\vert ]\nonumber\\
  & & + [\vert f_2^{(4)}(x+\theta _6 h)\vert +
      4 \vert \omega \vert \vert f_1^{(3)}(x)\vert ]
      \}\,.
\label{m2h2}
\end{eqnarray}
\subsection{A numerical example}

We consider the case study provided by the test function~(2.30) of
%                     \cite Ixaru [1997]
\cite{Ixaru97},
namely,
\begin{equation}
  \Phi (x) = f(x) \cos (\omega x)\,,\quad
  f(x) = 1/(1 + x)\,.
\label{ydtest}
\end{equation}
and compute, similar to
%                     \cite Ixaru [1997]
\cite{Ixaru97},
the errors associated to the discretized Leibniz and exponentially fitted
approximate derivation formulas at $x = 1$, under a step-size $h = 0.1$.
The computation is done for frequency parameter values $\omega \in [0, 80]$,
using the uniform sampling
\begin{equation}
  \omega _k = k h_{\omega}\,,\quad
  h_{\omega} = 0.1\,.
\label{osampling}
\end{equation}

In figures~1 and~2 of
%                     \cite Ixaru [1997]
\cite{Ixaru97},
the {\sl linearly scaled errors}
\begin{equation}
   \Delta ^{\prime}_{d}  =
     \left \{
       \begin{array}{ll}
          \Phi'(1)-\Phi_{d}'(1)          & \mbox{$\vert \omega \vert \le 1$}\\
          (\Phi'(1)-\Phi_{d}'(1))/\omega & \mbox{otherwise}, \quad d=2,4\,,
       \end{array}
     \right .
\label{scaleder1}
\end{equation}
have been plotted. To allow straightforward comparison of the two methods, the
$\omega$-dependence of the scaled errors~(\ref{scaleder1}) are shown in
figures~\ref{fig:ypoh2} and~\ref{fig:ypoh4} below.

To demonstrate the {\sl uniform bounds}~(\ref{mph2}) and~(\ref{mph4})
respectively, the absolute errors associated to the
${\cal O}(h^d)$ ($d\! =\! 2,\, 4$) Leibniz formulas~(\ref{ypoh2})
and~(\ref{ypoh4}),
\begin{equation}
   E^{\prime}_{d} = \Phi'(1)-\Phi_{d}'(1)\,,
\label{e1d}
\end{equation}
are also plotted in figures~\ref{fig:ypoh2} and~\ref{fig:ypoh4}.

\begin{figure}[h]
%\resizebox{\textwidth}{!}
\resizebox{\textwidth}{10cm}
{\includegraphics{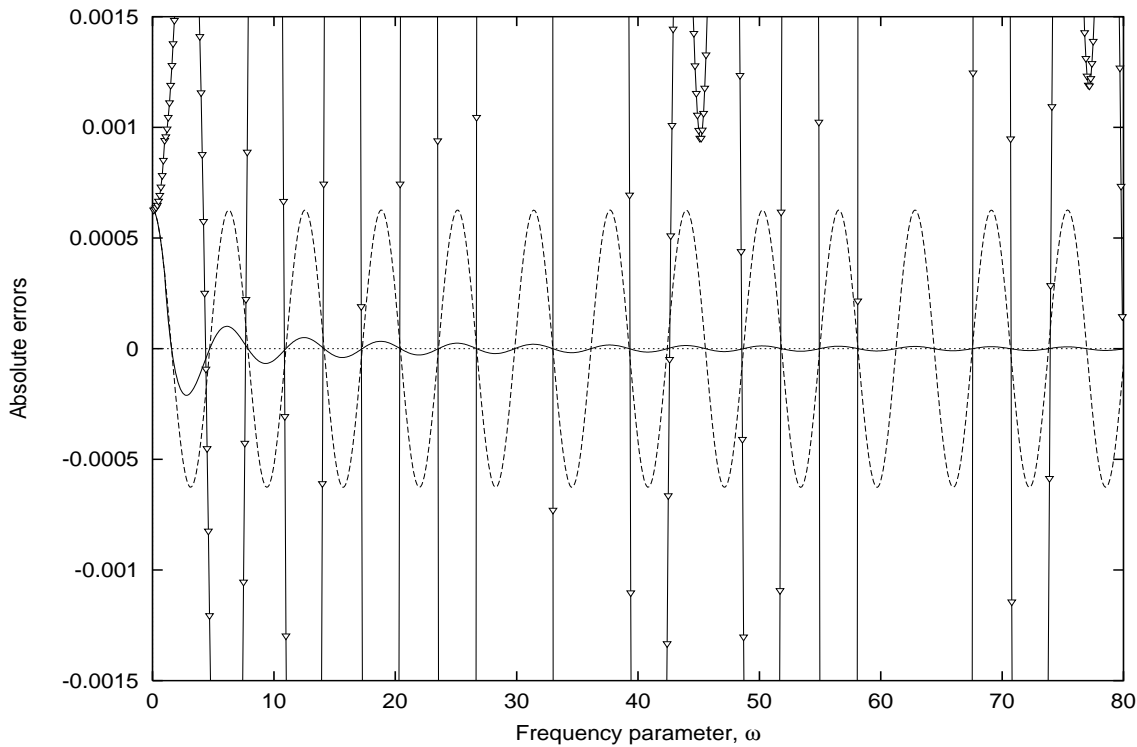}}
\caption{The $\omega$ dependence of the accuracies of the Leibniz derivation
         formula~(\ref{ypoh2}) and of the
         two-point exponentially fitted formula (2.28) of~%
%                     \cite Ixaru [1997]
\protect\cite{Ixaru97}.
         The linearly scaled absolute errors of~(\ref{ypoh2}) (solid line)
         damp out as $\omega$ increases. The
         linearly scaled errors of the exponentially fitted formula
         (linespoints) are small over narrow ranges of $\omega$ values only.
         The genuine absolute errors~(\ref{e1d}) of the Leibniz
         derivation formula~(\ref{ypoh2}) (dashed line) demonstrate the
         uniform bound~(\ref{mph2}).
        }
\label{fig:ypoh2}
\end{figure}

\begin{figure}[ht]
%\resizebox{\textwidth}{!}
\resizebox{\textwidth}{10cm}
{\includegraphics{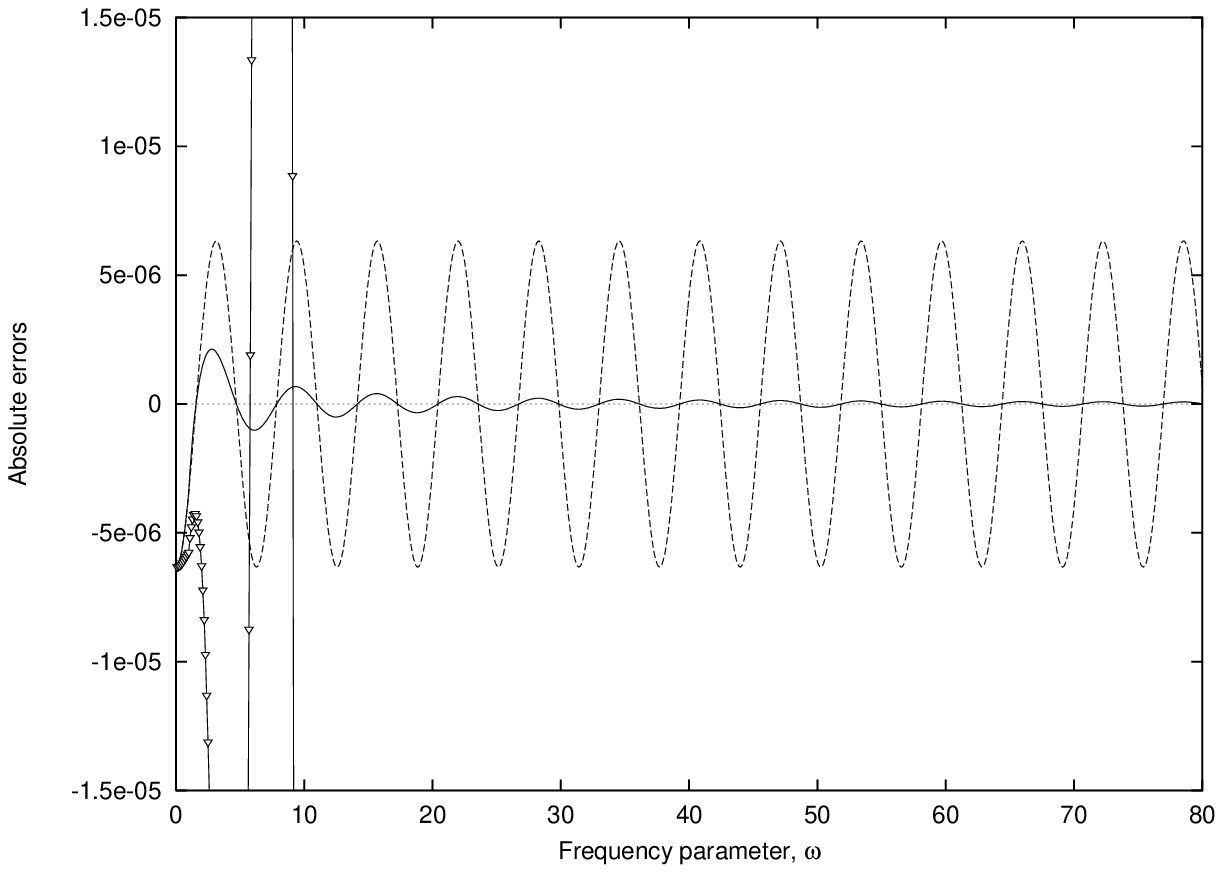}}
\caption{The $\omega$ dependence of the accuracies of the Leibniz derivation
         formula~(\ref{ypoh4}) and of the
         four-point exponentially fitted formula (3.1), (3.3) of~%
%                     \cite Ixaru [1997]
\protect\cite{Ixaru97}.
         The linearly scaled absolute errors of~(\ref{ypoh4}) (solid line)
         damp out as $\omega$ increases. Over the sampling~(\ref{osampling}),
         the linearly scaled errors of the exponentially fitted formula
         (linespoints) fall within the specified $y$-range (and are hence
         visible) at low $\omega$ values
         only. The genuine absolute errors~(\ref{e1d}) of the Leibniz
         derivation formula~(\ref{ypoh4}) (dashed line) demonstrate the
         uniform bound~(\ref{mph4}).
        }
\label{fig:ypoh4}
\end{figure}

There are three characteristic features which follow from the obtained data.
\begin{itemize}
  \item{The absolute errors associated to the Leibniz derivation formulas~%
        (\ref{ypoh2}) and (\ref{ypoh4}) are {\sl finite everywhere},
        irrespective of the value of the frequency parameter $\omega $, in
        agreement with the bounds~(\ref{mph2}) and~(\ref{mph4}).

        The magnitudes of the absolute errors associated to the exponentially
        fitted derivation formulas of reference~%
%                     \cite Ixaru [1997]
        \cite{Ixaru97}
        get {\sl arbitrarily large\/} in the neighbourhood of the critical
        points $\omega \! = \! k\pi \,,~k \! = \! \pm 1,~\pm 2,~\ldots $
       }
  \item{Let
        \begin{equation}
          \omega _k = \omega _0 + 2 k\pi \,,~k \! = \! 0,~\pm 1,~\pm 2,~\cdots 
          \label{om2k}
        \end{equation}
        For each of the two Leibniz approximate derivation formulas, Eqs.~%
        (\ref{ypoh2}) and (\ref{ypoh4}) respectively, Eqs.~(\ref{erh2})
        and~(\ref{erh4}) yield
        \begin{equation}
          E^{\prime}_d(\omega _k) = E^{\prime}_d(\omega _0)\,,
          \label{eperiod2k}
        \end{equation}
        i. e., the absolute errors~(\ref{e1d}) {\sl remain the same\/} over
        the set of equally spaced values~(\ref{om2k}). These errors show
        therefore {\sl periodic behaviours}, of {\sl constant specific
        amplitudes}, with respect to the frequency parameter $\omega $.
        Thus, in Fig.~\ref{fig:ypoh2}, the amplitude $0.627 \times 10^{-3}$
        equals the accuracy of the two-point derivation formula of the factor
        $f(x)$ in the test function~(\ref{ydtest}), while in
        Fig.~\ref{fig:ypoh4}, the amplitude
        $0.633 \times 10^{-5}$ equals the accuracy of the four-point
        derivation formula of the same factor $f(x)$.

        By contrast, the magnitudes of the absolute errors associated to the
        exponentially fitted formulas of the first order derivatives
        {\sl linearly increase with $\vert \omega \vert $\/} over the
        set~(\ref{om2k}). Hence the accuracy of these exponentially fitted
        formulas {\sl linearly deteriorates\/} with the increase of 
        $\vert \omega \vert $.
       }
  \item{If in Eqs.~(\ref{erh2}) and~(\ref{erh4}) the values $\theta _{1} =
        \theta _{2} = \theta _{3} = \theta _{4} = 1$ are assumed, then
        {\sl leading error estimates\/} are obtained. At the scale of the
        present figures~\ref{fig:ypoh2} and~\ref{fig:ypoh4}, the exact
        errors~(\ref{e1d}) are practically {\sl indistinguishable from
        these estimates\/} (compare the amplitudes
        $A_2(exact) = 0.627 \times 10^{-3}$ with
        $A_2(est.) = 0.625 \times 10^{-3}$, and
%       $A_2(estimated) = 0.625 \times 10^{-3}$, and
        $A_4(exact) = 0.633 \times 10^{-5}$ with
        $A_4(est.) = 0.625 \times 10^{-5}$). This is a consequence of
%       $A_4(estimated) = 0.625 \times 10^{-5}$). This is a consequence of
        the particular test function (\ref{ydtest}) proposed in~%
%                     \cite Ixaru [1997]
        \cite{Ixaru97}
        and it should not be overemphasized.

        Figure~1 of~%
%                     \cite Ixaru [1997]
        \cite{Ixaru97},
        however, shows that, for the same case study, the estimated errors
        of the exponentially fitted formulas {\sl may vary significantly\/}
        from the exact error values.
       }
\end{itemize}

The bound (\ref{m2h2}) shows that the appropriate quantities for the study
of the dependence of the accuracy of the discretization errors~(\ref{ysoh2})
versus $\omega$ are the {\sl linearly scaled errors}
\begin{equation}
 E^{\prime \prime}_{2}  =
   \left \{
     \begin{array}{ll}
        \Phi^{\prime\prime}(1)-\Phi_{2}^{\prime\prime}(1) &
               \mbox{$\vert \omega \vert \le 1$}\\
        (\Phi^{\prime\prime}(1)-\Phi_{2}^{\prime\prime}(1))/\omega &
               \mbox{otherwise}
     \end{array}
   \right .
\label{lscaleder2}
\end{equation}

In figure~3 of 
%                     \cite Ixaru [1997]
\cite{Ixaru97},
however, the dependence versus $\omega$ of the accuracy of the exponentially
fitted three-point second order derivative has been plotted in terms of the
{\sl quadratically scaled errors}
\begin{equation}
 \Delta ^{\prime \prime}_{2}  =
   \left \{
     \begin{array}{ll}
        \Phi^{\prime\prime}(1)-\Phi_{2}^{\prime\prime}(1) &
               \mbox{$\vert \omega \vert \le 1$}\\
        (\Phi^{\prime\prime}(1)-\Phi_{2}^{\prime\prime}(1))/\omega ^{2} &
               \mbox{otherwise}
     \end{array}
   \right .
\label{qscaleder2}
\end{equation}

\begin{figure}[h]
%\resizebox{\textwidth}{!}
\resizebox{\textwidth}{10cm}
{\includegraphics{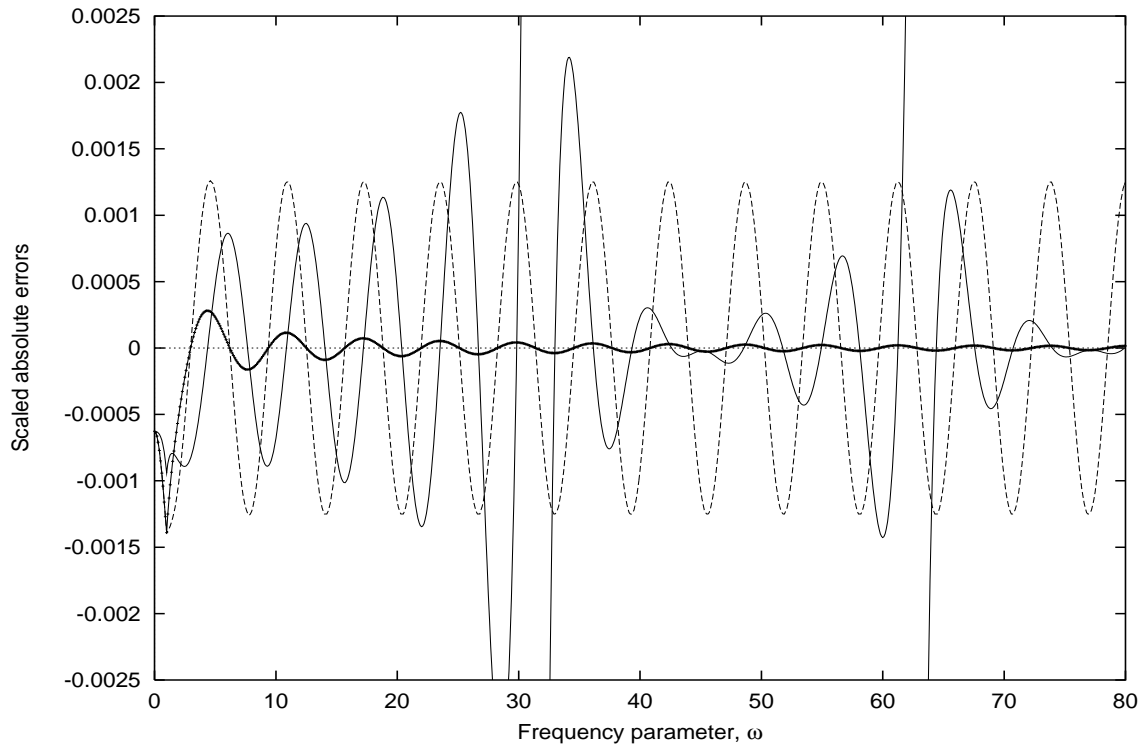}}
\caption{The $\omega$ dependence of the accuracies of the Leibniz derivation
         formula~(\ref{ysoh2}) and of the
         three-point exponentially fitted formula (3.16), (3.18) of reference~%
%                     \cite Ixaru [1997]
\protect\cite{Ixaru97}.
         The quadratically scaled absolute errors of~(\ref{ysoh2})
         (linespoints) damp out as $\omega$ increases. The
         quadratically scaled errors of the exponentially fitted formula
         (solid line) are shown for comparison. The linearly scaled
         absolute errors~(\ref{lscaleder2}) of the Leibniz
         derivation formula~(\ref{ysoh2}) (dashed line) show constant
         amplitude with $\omega$ at $\omega > 1$.
        }
\label{fig:ysoh2}
\end{figure}

Figure~\ref{fig:ysoh2} summarizes the $\omega$ dependence of the quadratically
scaled errors of the exponentially fitted and discrete Leibniz derivatives
as well as the linearly scaled errors of the discrete Leibniz derivative~%
(\ref{ysoh2}).

From these data, we draw the following conclusions on the approximation
of the second order derivatives within the two approaches.
\begin{itemize}
  \item{
        Similar to the result obtained for the first order derivatives, the
        errors associated to the Leibniz formula are {\sl finite everywhere},
        while those associated to the exponentially fitted formula
        {\sl diverge\/} at the special point set
        $\omega \! = \! k\pi \,,~k \! = \! \pm 1,~\pm 2,~\ldots $
       }
  \item{
        At frequency parameter values $\vert \omega \vert > 1$, the Leibniz
        formula is almost everywhere better than the exponentially fitted
        one. The set of arguments $\omega x + \delta$ at which the accuracy
        of the latter formula is comparable to that of the former gets
        monotonically smaller as the magnitude of $\omega $ increases, tending
        to a countable manifold at asymptotically large $\vert \omega \vert$.
       }
  \item{
        For the present case study, a phase difference, roughly equal to
        $\pi / 2$, is noticeable between the $\omega$-dependences of the
        errors associated to the two approaches. This feature stems from the
        prevalence of different asymptotic terms in the discretization errors:
        a term proportional to $\omega ^{2} \cos (\omega)$ in the
        exponentially fitted formula of~%
%                     \cite Ixaru [1997]
        \cite{Ixaru97},
        and a term proportional to $\omega \sin (\omega)$ in the Leibniz
        formula.
       }
\end{itemize}

The evidence accumulated for the first and second order derivatives raises
the question on the rate of deterioration of the accuracy of the
discretization of the $n$-th order derivative of the function~(\ref{reffunc})
under asymptotically large $\omega$ values.

        Within the exponential fitting, this rate gets proportional
        to $\vert \omega \vert ^n$, therefore this approach yields
        {\sl unsatisfactory results for derivatives of any order $n$}, even
        $n=1$, which corresponds to the first order derivatives.

        Let $f^{(m)}_{d}(x)$ denote the ${\cal O}(h^{d})$ difference formula
        of the $m$-th order derivative of the regular factor $f(x)$
        entering the reference function (\ref{ygeneric}). Then the
        discretization error associated to the ${\cal O}(h^{d})$
        approximation of equation~(\ref{Leibniz}) reads
        \begin{equation}
           E^{(n)}_{d}(x) = F^{(n)}(x) - F^{(n)}_{d}(x) = 
              \sum_{k=0}^{n-1}
                   \left(
                       \begin{array}{cc}
                         n\\
                         k
                        \end{array}
                   \right)
                     \omega^{k}[f^{(n-k)}(x) - f^{(n-k)}_{d}(x)]
                         g_{s,\eta}^{(k)}(t)\,,
        \label{err-d-Leibniz}
        \end{equation}
        therefore it behaves like $\vert \omega \vert ^{n-1}$ at
        asymptotically large $\omega $.
\section{Quadratures}
\label{sec:quad}

\subsection{Interpolatory Simpson quadrature sums\\ for oscillatory functions}
\label{sec:simpson}

Let
\begin{equation}
  I^{\omega,\delta}_{s,\eta}[a,b]f = 
               \int_a^b f(x) g_{s,\eta}(\omega x + \delta) d\!x,
               \; s = 1, 2,\; \eta = \pm 1\, ,
\label{qgeneric}
\end{equation}
denote the distinct integrals corresponding to each of the four basic
integrands~(\ref{ygeneric}).

In this section, we discuss interpolatory Simpson quadrature sums for
oscillatory functions. These provide an alternative to the exponentially
fitted formulas derived in section~4 of
%                     \cite Ixaru [1997]
\cite{Ixaru97}.
We start with the derivation of reliable quadrature sums for the solution of
the basic integrals~(\ref{qgeneric}). Simpson quadrature sums for integrals
\begin{equation}
  I[a,b]\Phi = \int_a^b\Phi(x) d\!x,
\label{iofy}
\end{equation}
with integrands $\Phi(x)$ of the form (\ref{reffunc}), will then be given by
linear combinations of the basic quadrature sums.

An important preliminary operation is the mapping of the integration range
$[a, b]$ onto the fundamental range $[-1, 1]$ via the transform of
independent variable,
\begin{equation}
  x = c + h y; \; y\in [-1, 1]\, ;\; c = (b + a)/2\, ;\; h = (b - a)/2\, ,
\label{chvar}
\end{equation}
where $c$ and $h$ denote the {\sl centre\/} and
respectively the {\sl half-length\/} of the integration range $[a, b]$.

As a result of this operation, we get the following relationships which
are best written in matrix form
\begin{equation}
 \left(
   \begin{array}{c}
     I^{\omega,\delta}_{1,\eta}[a,b]f\\
     I^{\omega,\delta}_{2,\eta}[a,b]f\\
   \end{array}
 \right) = 
 \left(
   \begin{array}{cc}
     R_{1,\eta}&\eta R_{2,\eta}\\
     R_{2,\eta}&R_{1,\eta}\\
    \end{array}
 \right)
 \left(
    \begin{array}{c}
      I^\lambda_{1,\eta}[-1,1]\psi \\
      I^\lambda_{2,\eta}[-1,1]\psi \\
     \end{array}
 \right)
        ; \quad \eta = \pm 1\, .
\label{iab-11}
\end{equation}
Thus, the calculation of each of the two integrals~(\ref{qgeneric}) occurring
in the left hand side has been reduced to the evaluation of the {\sl pair\/}
of fundamental integrals
\begin{equation}
  I^{\lambda}_{s,\eta}[a,b]\psi = 
               \int_{-1}^{1} \psi (y) g_{s,\eta}(\lambda y) d\!y, \quad 
               \psi (y) = f(c + hy)\,; \quad s = 1, 2\, ,
\label{qfundam}
\end{equation}
with a transfer matrix R of elements
\begin{equation}
  R_{s,\eta}=h g_{s,\eta}(\varphi _{c})\,,\quad s=1,2\,;\; \eta=\pm1 \,.
\label{eqR}
\end{equation}
Here, $\lambda$ is the dimensionless parameter defined
by Eq.~(\ref{deflambda}), while $\varphi _c = \omega c + \delta$ denotes the
phase of the center $c$ of the integration range $[a, b]$.

To derive Simpson quadrature sums for each of the basic
integrals~(\ref{qgeneric}), the regular part $f(x)$ of the
product~(\ref{ygeneric}) is {\sl interpolated\/} by an arc of parabola,
$L_2(x)$, at the knots
\begin{equation}
  a = x_0\ <\ x_2\ <\ x_1 = b\, .
\label{qknots}
\end{equation}
We use the Newton representation of the unique Lagrange interpolation polynomial
of $f(x)$ at the knots (\ref{qknots})
%                     \cite deBoor [1978]
\cite{deBoor78}.

The generalization of the classical Simpson quadrature to the
solution of the integral (\ref{qgeneric}) will then be
\begin{equation}
  Q^{\omega,\delta}_{s,\eta}[a,b]f \equiv
  I^{\omega,\delta}_{s,\eta}[a,b]L_2 = \int_a^b L_2(x)
     g_{s,\eta}(\omega x + \delta) d\!x \, .
\label{ql2x}
\end{equation}

Taking into account~(\ref{iab-11}), the computation of the
last integral is reduced to that of the pair of integrals
$I^\lambda_{1,\eta}[-1,1]\ell_2$ and $I^\lambda_{2,\eta}[-1,1]\ell_2$,
where $\ell_2(y) = L_2(c + hy)$. To get well-conditioned expressions of these
integrals, $\ell_2(y)$ is expressed in terms of Chebyshev polynomials of the
first kind,
\begin{equation}
  \ell_2(y) = \sum_{i=0}^{2} \beta_{i2} T_{i}(y)\, ,
\label{l2y}
\end{equation}
with the coefficients given respectively by
\begin{eqnarray}
  \beta_{02} &\!\! = \!\!& {\scriptstyle{\frac{1}{8}}}[
          (3 - 1/\rho) f_0 +
          (2 + \rho + 1/\rho) f_2 +
          (3 - \rho) f_1]\, ,\nonumber\\
  \beta_{12} &\!\! = \!\!& {\scriptstyle{\frac{1}{2}}}(f_1 - f_0)\, ,\nonumber\\
  \beta_{22} &\!\! = \!\!& {\scriptstyle{\frac{1}{8}}}[
          (1 + 1/\rho) f_0 -
          (2 + \rho + 1/\rho) f_2 +
          (1 + \rho) f_1]\, .
\label{cl2y}
\end{eqnarray}
Here
\begin{equation}
  f_k \equiv f(x_k),\; k = 0, 1, 2\\
\label{findicek}
\end{equation}
while
\begin{equation}
  \rho = \frac{1+y_2}{1-y_2},
\label{defrho}
\end{equation}
denotes the ratio of the lengths of the two subranges created by the
inner abscissa $y_2 = (x_2 - c)/h$ inside the fundamental range [-1, 1].

General expressions of the reduced integrals
$I_{s,\eta}^{\lambda}[-1,1]\ell_n$ have been derived in~%
%                     \cite Adam [1991]
\cite{Adam91}.
[There, however, two misprints have to be corrected. First, the right hand side
of the expression (28) of the coefficients $s_{q,2k-1}$ is to be preceded
by a minus sign. Second, in Eq.~(39), the summation is to proceed from 1 to $k$,
over hypergeometric functions
$_0F_1(j+{\scriptstyle{\frac{3}{2}}};
{\scriptstyle{\frac{1}{4}}}\eta \lambda^2)$.]
In view of the low polynomial degree of the expansion (\ref{l2y}), equations
(38) and (39) of
%                     \cite Adam [1991]
\cite{Adam91}
provide suitable expressions for the moments of the involved Chebyshev
polynomials.

The obtained expression of $I^{\lambda}_{1,\eta}[-1,1]\ell_2$ is transformed
to a well-conditioned one by use of the recurrence relation (40) of
%                     \cite Adam [1991]
\cite{Adam91}
and we finally get
\begin{eqnarray}
  I^\lambda_{1,\eta}[-1,1]\ell_{2} &\!\! = \!\!&
              C_1(\rho )\cdot {_0F_1}({\scriptstyle{\frac{5}{2}}};
        {\scriptstyle{\frac{1}{4}}}\eta \lambda ^{2}) +
   \frac{1}{15}\eta \lambda ^{2}(f_0 + f_1)
        {_0F_1}({\scriptstyle{\frac{7}{2}}};
        {\scriptstyle{\frac{1}{4}}}\eta \lambda ^{2})\, ,
\label{redi1}\\
  I^\lambda_{2,\eta}[-1,1]\ell_{2} & \!\!= \!\!&
   \frac{1}{3}\lambda (f_1 - f_0)
        {_0F_1}({\scriptstyle{\frac{5}{2}}};
        {\scriptstyle{\frac{1}{4}}}\eta \lambda ^{2})\, ,
\label{redi2}
\end{eqnarray}
where
\begin{equation}
   C_1(\rho ) = \frac{1}{3}\Bigl[(2 - 1/{\rho}) f_0 +
                (2 + \rho + 1/{\rho}) f_2 + (2 - \rho) f_1\Bigr]\,.
\label{c1rho}
\end{equation}
If the knots~(\ref{qknots}) are {\sl equally spaced}, then
$\rho = 1$ and the coefficient~(\ref{c1rho}) goes into the classical Simpson
quadrature sum for the regular factor $f(x)$,
\begin{equation}
  C_1(1) = \frac{1}{3}\Bigl(f_0 + 4 f_1 + f_2\Bigr)\,.
\label{c11}
\end{equation}

Similar to the interpolatory quadrature sums at Clenshaw-Curtis and related
abscissas,
%                     \cite Adam [1991]
\cite{Adam91},
the only practical difficulty in the implementation of these quadrature sums
into a code is the accurate computation of the hypergeometric functions
$_0F_1$ of interest. The preparation of a paper which describes the computer
programme devoted to the computation of basis sets of hypergeometric
functions~$_0F_1$ to machine accuracy is underway and it will be submitted
to this journal in the nearest future.

Equations~(\ref{ql2x}), (\ref{iab-11}), (\ref{redi1}), (\ref{redi2}),
(\ref{c11}) provide the required generalization of the classical Simpson
quadrature rule to each of the four reference functions~(\ref{ygeneric}).
Indeed, in the limit $\omega \rightarrow 0$, these equations result in
$Q^{0,\delta}_{s,\eta}[a,b]f = 
{\scriptstyle{\frac{1}{3}}}h(f_0 + 4 f_2 + f_1)g_{s,\eta}(\delta)$. On the
other side, in the same limit, Eq.~(\ref{ygeneric}) yields
$F(x) = f(x) g_{s,\eta}(\delta)$, hence $F_k = f_k g_{s,\eta}(\delta)$.
Therefore, $Q^{0,\delta}_{s,\eta}[a,b]f =
{\scriptstyle{\frac{1}{3}}}h(F_0 + 4 F_2 + F_1)$, which is nothing but
the classical Simpson quadrature sum $Q[a, b]F$.

The Simpson quadrature sum for the integral (\ref{iofy})
with the integrand $\Phi(x)$ given by (\ref{reffunc}) is then given by
\begin{eqnarray}
  Q^{\omega,\delta}_{\eta}[a,b]\Phi
     &\!\! = \!\!& \Bigl(I^\lambda_{1,\eta}[-1,1]\ell_{2}^{\prime} +
      I^\lambda_{2,\eta}[-1,1]\ell_{2}^{\prime\prime}\Bigr)R_{1,\eta}\nonumber\\
      &\!\! + \!\!& \Bigl(I^\lambda_{1,\eta}[-1,1]\ell_{2}^{\prime\prime} +
      \eta I^\lambda_{2,\eta}[-1,1]\ell_{2}^{\prime}\Bigr)R_{2,\eta}\, ,
\label{qofy}
\end{eqnarray}
where $\ell_{2}^{\prime}(y)$ and $\ell_{2}^{\prime\prime}(y)$ denote
the second degree interpolatory polynomials associated to the
regular factors $f_1(c+hy)$ and $f_2(c+hy)$ respectively.
\subsection{Leading error estimate}
\label{sec:qerror}

The error associated to the quadrature sum~(\ref{ql2x}), which replaces
the computation of the basic integral~(\ref{qgeneric}), is given by
\begin{eqnarray}
  I^{\omega,\delta}_{s,\eta}[a,b]\Delta f &\!\! = \!\!& 
        \int_a^b \Delta f(x) g_{s,\eta}(\omega x + \delta) d\!x,\nonumber\\
  \Delta f(x) &\!\! = \!\!& f(x) - L_2(x)\,.
\label{qrest}
\end{eqnarray}

In view of the relationship~(\ref{iab-11}), an estimate of this expression
requires two estimates~(\ref{qfundam}), $s = 1, 2$, for the function
\begin{equation}
  \Delta \psi(y) = \psi(y) - \ell_2(y) = f(c + hy) - \ell_2(y).
\label{deltapsi}
\end{equation}

To get the leading terms of the two estimates, the reference function
$f(c + hy)$ is expanded in Taylor series up to fourth order around the point
$x = c$ while the coefficients of the interpolation polynomial are expressed
in terms of the ${\cal O}(h^{2})$ approximate derivatives $f^{\prime}_{2}(c)$
and $f^{\prime\prime}_{2}(c)$ of $f(x)$. After straightforward algebra, we get
the Chebyshev polynomial expansion
\begin{equation}
  \Delta \psi(y) \approx \frac{1}{2}\alpha_{0}[T_{2}(y)-T_{0}(y)] +
                         \frac{1}{8}\alpha_{4}[T_{4}(y)-T_{0}(y)] +
                         \frac{1}{4}\alpha_{3}[T_{3}(y)-T_{1}(y)]\,,
\label{apprdeltapsi}
\end{equation}
where
\begin{eqnarray}
  \alpha_{0} &\!\! = \!\!& \frac{1}{4}[ (1 - 1/\rho) f_0 +
          (2 + \rho + 1/\rho) f_2 + (1 - \rho) f_1] - f(c)\nonumber\\
             &\!\! \approx \!\!& - \frac{1}{4}(3 + \rho )y_2 \alpha_{3} 
          - y_2^2 \alpha_{4}\,,
\label{defalfa0}\\
  \alpha_{3} &\!\! = \!\!& \frac{1}{6} h^3 f^{(3)}(c)\,,
\label{defalfa3}\\
  \alpha_{4} &\!\! = \!\!& \frac{1}{24} h^4 f^{(4)}(c)\,.
\label{defalfa4}
\end{eqnarray}
The coefficient $\alpha_0$ vanishes identically in the case of equally
spaced knots. It arises only within the generalized Simpson quadrature sum~%
(\ref{redi1}) under $\rho \ne 1$.

The leading error estimates associated to the quadrature sums~(\ref{redi1})
and~(\ref{redi2}) are now immediate:
\begin{eqnarray}
  I^\lambda_{1,\eta}[-1,1]\Delta \psi &\!\! = \!\!&
        -\frac{4}{3}\Bigl[
        (\alpha _0 + \alpha _4) \cdot {_0F_1}({\scriptstyle{\frac{5}{2}}};
        {\scriptstyle{\frac{1}{4}}}\eta \lambda ^{2}) -
   \frac{4}{5} \alpha _4 \cdot
        {_0F_1}({\scriptstyle{\frac{7}{2}}};
        {\scriptstyle{\frac{1}{4}}}\eta \lambda ^{2})
                    \Bigr] \, ,
\label{eri1}\\
  I^\lambda_{2,\eta}[-1,1]\Delta \psi &\!\! = \!\!&
   -\frac{4}{15}\lambda \alpha _3 \cdot
        {_0F_1}({\scriptstyle{\frac{7}{2}}};
        {\scriptstyle{\frac{1}{4}}}\eta \lambda ^{2})\, .
\label{eri2}
\end{eqnarray}

Taking into account the expressions~(\ref{defalfa0}),~(\ref{defalfa3}), and~%
(\ref{defalfa4}), we conclude that the equally spaced Simpson quadrature sums~%
(\ref{redi1}) and~(\ref{redi2}) show ${\cal O} (h^4)$ accuracy. However,
under a nonuniform mesh, the order of accuracy is reduced to ${\cal O} (h^3)$,
due to the non-vanishing $\alpha_0$ coefficient.

Let in~(\ref{eri1}) and~(\ref{eri2}) the hypergeometric functions
${_0F_1}({\scriptstyle{\frac{5}{2}}};
{\scriptstyle{\frac{1}{4}}}\eta \lambda ^{2})$ and
${_0F_1}({\scriptstyle{\frac{7}{2}}};
{\scriptstyle{\frac{1}{4}}}\eta \lambda ^{2})$
be expressed in terms of the fundamental hypergeometric functions
${_0F_1}({\scriptstyle{\frac{1}{2}}};
{\scriptstyle{\frac{1}{4}}}\eta \lambda ^{2})$ (which equals $\cos(\lambda)$
under $\eta = - 1$, respectively $\cosh(\lambda)$ under $\eta = 1$) and
${_0F_1}({\scriptstyle{\frac{3}{2}}};
{\scriptstyle{\frac{1}{4}}}\eta \lambda ^{2})$ (which equals
$\sin(\lambda) / \lambda$ under $\eta = - 1$, respectively
$\sinh(\lambda) / \lambda$ under $\eta = 1$).
These operations put into evidence the occurrence of an {\sl $\omega ^{-2}$
power law at asymptotically large $\omega$}, coming from the hypergeometric
functions. Therefore, the accuracy of the interpolatory Simpson quadrature
sums is expected to {\sl increase\/} with increasing $\omega$.
\subsection{An accuracy test}

A meaningful comparison of the above interpolatory Simpson and exponentially
fitted quadrature sums of~%
%                     \cite Ixaru [1997]
\cite{Ixaru97}
is provided by the test case (4.18) of reference
%                     \cite Ixaru [1997]
\cite{Ixaru97},
that is, the integral~(\ref{iofy}) of the oscillatory function
\begin{equation}
  \Phi(x)=-f^2(x)\cos (\omega x)-\omega f(x)\sin (\omega x),\ f(x)=1/(1+x)\, ,
\label{quadfunc}
\end{equation}
over the range $[a, b] = [0.9, 1.1]$.
The primitive of the integrand (\ref{quadfunc}) is the function~(\ref{ydtest}),
hence the exact value of $I[a, b]\Phi$, Eq.~(\ref{iofy}), is immediate:
\begin{equation}
  I[a, b]\Phi = - \Bigl[\frac{h}{2}\cos(\varphi_{c}) \cos(\lambda) +
                        \frac{1+c}{2}\sin(\varphi_{c}) \sin(\lambda)
                  \Bigr]/\Bigl[\left(\frac{1+c}{2}\right)^{2} -
                         \left(\frac{h}{2}\right)^{2}\Bigr]\,.
\label{exacti}
\end{equation}
where $c\! =\! 1\,;\ h\! =\! 0.1\,.$
With the quadrature sum $Q[a, b]\Phi$ given by~(\ref{qofy}), we can
compute the absolute quadrature errors, at various values of the frequency
parameter $\omega$,
\begin{equation}
  \Delta = I[a, b]\Phi-Q[a, b]\Phi \,.
\label{erq3p}
\end{equation}

\begin{figure}[h]
%\resizebox{\textwidth}{!}
\resizebox{\textwidth}{10cm}
{\includegraphics{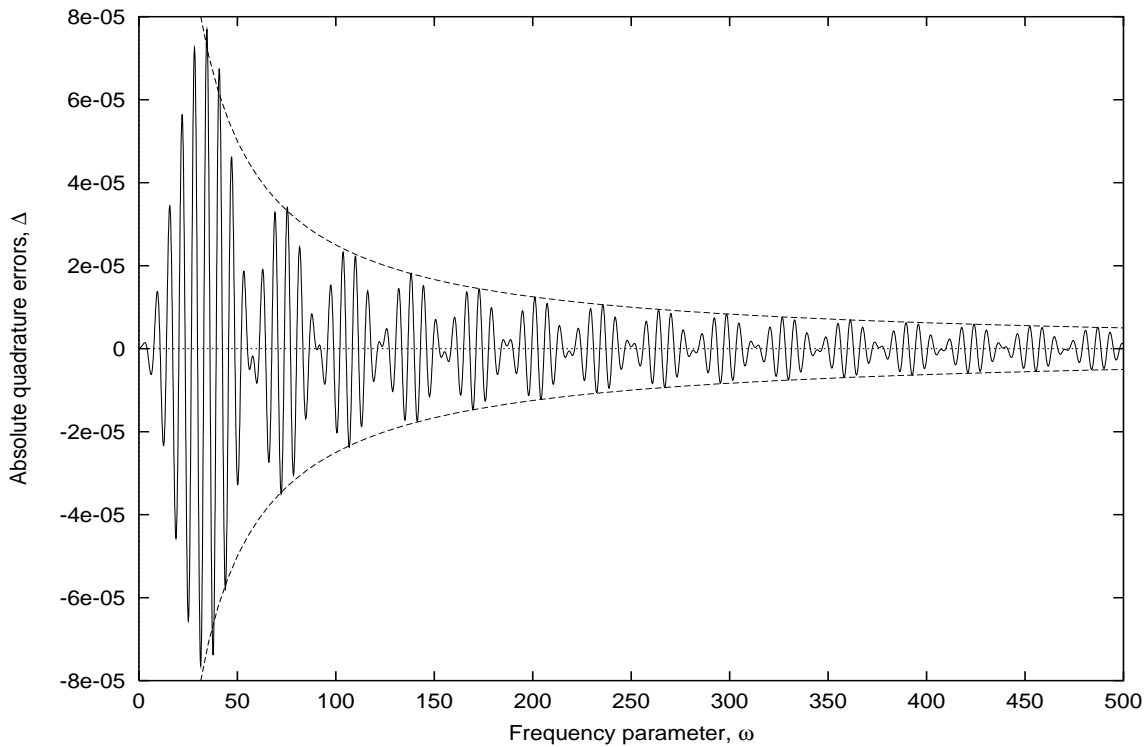}}
\caption{Absolute quadrature errors, $\Delta$, Eq.~(\ref{erq3p}), associated
         to the interpolatory Simpson method for the oscillatory function~%
         (\ref{quadfunc}) (solid line). The two envelopatrices
         $\Delta = \pm 0.0025 / \omega$ (dashed lines) illustrate the
         asymptotic error behaviour established by the error analysis.
        }
\label{fig:isimpson}
\end{figure}

An important feature stemming from equation~(\ref{exacti}) is the occurrence
of two distinct periodic oscillatory behaviours of the considered integral.
First, there is a characteristic {\sl $\varphi _c$-dependent periodicity},
defined by the phase $\varphi _c = \omega c + \delta$ of the centre $c$ of
the integration range through the factors $\cos(\varphi_{c})$ and
$\sin(\varphi_{c})$. Here, this kind of periodicity is characterized by
{\sl short wavelength\/} oscillations of period $T_{\varphi_{c}} = 2 \pi$.
Second, there is a characteristic {\sl $\lambda$-dependent periodicity},
defined by the phase $\lambda$ of the integration range half-length through
the factors $\cos(\lambda )$ and $\sin(\lambda )$. Here, this kind of
periodicity is characterized by {\sl long wavelength\/} oscillations, of
period $T_{\lambda} = 2 \pi / h = 20 \pi$, which is ten times larger than
$T_{\varphi_{c}}$.

\begin{figure}[h]
%\resizebox{\textwidth}{!}
\resizebox{\textwidth}{10cm}
{\includegraphics{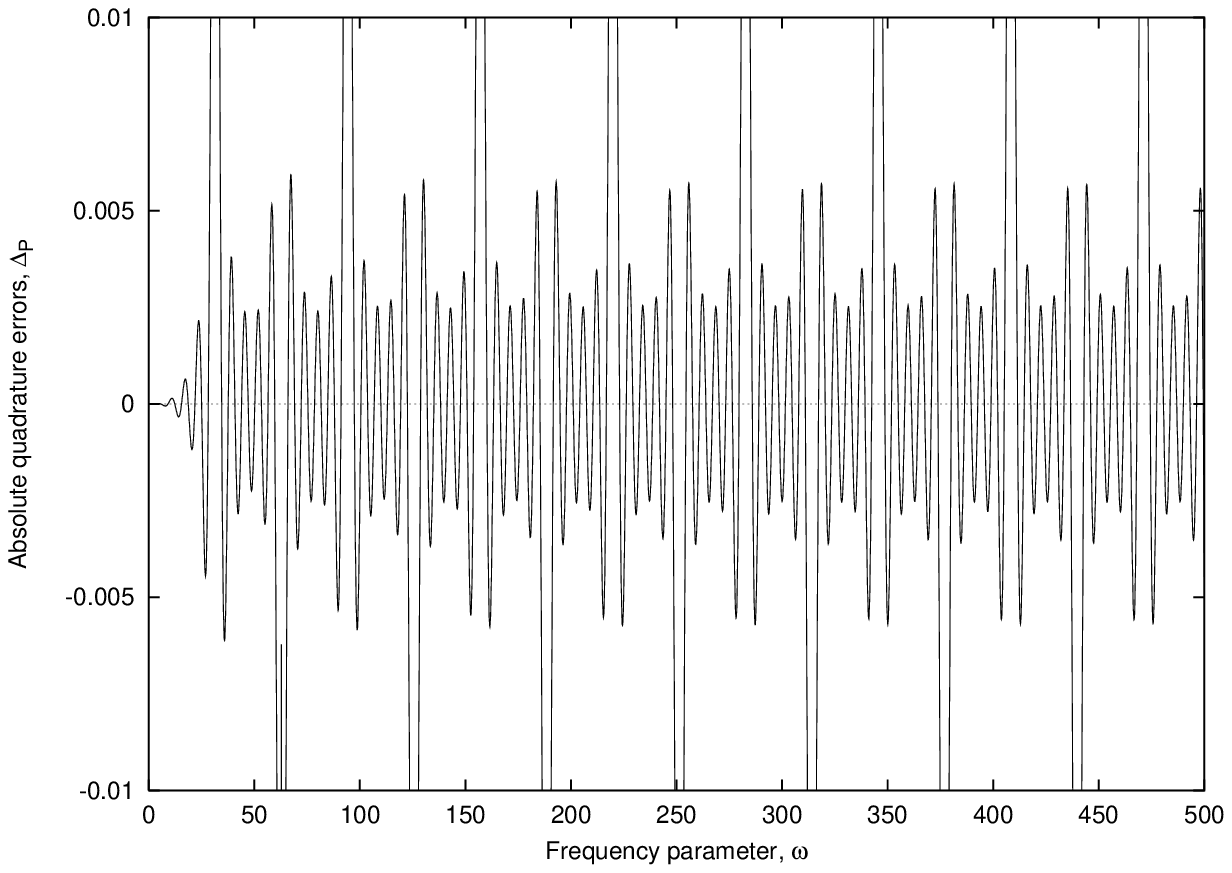}}
\caption{Absolute quadrature errors associated to the optimal exponentially
         fitted $P = 1$ method of~%
%                     \cite Ixaru [1997]
\protect\cite{Ixaru97}.
        }
\label{fig:isexfit}
\end{figure}

Figure~\ref{fig:isimpson} plots the obtained dependence $\Delta$ versus
$\omega$ over the frequency parameter range $\omega \in [0, 500]$.
To allow straightforward
comparison with the optimal exponentially fitted $P\!=\!1$ method of~%
%                     \cite Ixaru [1997]
\cite{Ixaru97},
a supplementary plot concerning the $\omega$ dependence of the quadrature
errors of this method is given in figure~\ref{fig:isexfit}. In both figures,
the sampling~(\ref{osampling}) was used. Separate plots have been needed
simply because the relevant $y$-scale of figure~\ref{fig:isimpson},
$(-8.0\times 10^{-5},\, +8.0\times 10^{-5})$, is two orders of magnitude
smaller than that required in figure~\ref{fig:isexfit},
$(-1.0\times 10^{-2},\, +1.0\times 10^{-2})$.

The inspection of the obtained numerical evidence results in the following
conclusions concerning the interpolatory Simpson and exponentially fitted
quadrature sums.
\begin{itemize}
  \item{Both kinds of periodicities discussed above are present in
        figures~\ref{fig:isimpson} and~\ref{fig:isexfit}. While the
        $\varphi _c$-induced periodicity pattern is the same in the two
        figures, the $\lambda$-induced periodicity pattern is substantially
        different.
       }
  \item{Being controlled by the function
        $\eta_0(-\lambda ^2) = \sin(\lambda ) /\lambda$, the $\lambda$-%
        periodicity of the exponentially fitted method shows, over sets of
        points $\omega \geq T_{\lambda}/2$, a period {\sl exactly\/} equal
        to $T_{\lambda}$.

        The error magnitudes are {\sl roughly the same\/} over sets of
        $\omega$ values separated by multiples of $T_{\lambda}$.
       }
  \item{By contrast, the $\lambda$-periodicity of the interpolatory Simpson
        quadrature sums is controlled by higher order hypergeometric
        functions $_0F_1$, the periods of which only {\sl asymptotically\/}
        tend towards the value $T_{\lambda}$. (Thus, for the present problem,
        approximately periodic long wavelength oscillations can be defined
        starting with values $\omega > T_{\lambda}$. As $\omega$
        increases, the quasi-period lengths {\sl monotonically decrease\/}
        towards $T_{\lambda}$. At values $\omega > 400$ an agreement with
        $T_{\lambda}$ within values better than one percent is obtained.)

        The behaviour of the error amplitudes with $\omega$ is in agreement
        with the error analysis done in subsection~\ref{sec:qerror}.
        The important detail to be noted here is the occurrence, in the
        integrand~(\ref{quadfunc}), of a non-oscillatory term which is
        proportional to $\omega$. This results into {\sl an asymptotic
        $\omega^{-1}$ power law behaviour\/} of the $\lambda$-defined long
        wavelength amplitudes. In figure~\ref{fig:isimpson}, the occurrence of
        this asymptotic regime is visualized by means of the two envelopatrices
        $\Delta = \pm 0.0025 / \omega$, which provide true upper bounds to the
        discretization errors at $\omega$ values beyond $T_{\lambda}$.
       }
  \item{In reference
%                     \cite Ixaru [1997]
        \cite{Ixaru97},
         breakdown of the quadrature errors was noticed to occur
         around $\lambda = 2k\pi$ for the suboptimally fitted method $P\!=\!0$
         and around the critical values $\lambda = k\pi$ $(k=1,\ 2,\ ...)$,
         for the optimally fitted method $P\!=\!1$ (only the latter being
         visible at the error scale, $-0.01 < \Delta_P < 0.01$, of fig.~4
         of reference
%                     \cite Ixaru [1997]
         \cite{Ixaru97}%
         ). Figure~\ref{fig:isexfit} is consistent with this frame.
         An increase of the mesh density~(\ref{osampling}) adds supplementary
         details in the neighbourhood of the critical points only.
       }
\end{itemize}
\section{Conclusions}
\label{sec:concl}

The approximate $p$-point Leibniz derivation formulas and the
interpolatory Simpson quadrature sums discussed in this paper provide simple,
versatile, and efficient tools for the analysis of oscillatory phenomena.

Both theoretical considerations and numerical evidence concerning
the dependence of the discretization errors on the frequency parameter of
the oscillatory functions show that the accuracy gain of the present formulas
over those based on the exponential fitting is overwhelming:
\begin{enumerate}
  \item{The discretization errors of the present numerical differentiation
        and integration formulas are {\sl finite everywhere\/} with respect
        to the frequency parameter $\omega$, whereas those of the
        corresponding exponentially fitted ones {\sl diverge\/} at specific
        countable sets of $\omega$ values.
       }
  \item{Even if the huge errors coming from convenient neighbourhoods around
        diverging points of the exponentially fitted methods are disregarded,
        the accuracies of the present methods still remain, in the average,
        {\sl two orders of magnitude better\/} with respect to those based
        on the exponential fitting.
       }
  \item{Over sets of $\omega$ values separated by entire periods of the
        oscillatory function, the absolute errors of the approximate Leibniz
        $n$-th order derivatives increase following an
        $\vert \omega \vert ^{n-1}$ power law, as compared to the
        $\vert \omega \vert ^{n}$ power law specific to the exponentially
        fitted formulas.

        Thus, the Leibniz first order derivatives~(\ref{ypoh2})
        and~(\ref{ypoh4}) are characterized by the {\sl uniform bounds\/}~%
        (\ref{mph2}) and~(\ref{mph4}) respectively, whereas the exponentially
        fitted ones {\sl linearly deteriorate\/} with $\omega$. As it concerns
        the second order derivative, the Leibniz formula~(\ref{ysoh2}) shows
        {\sl linear deterioration with~$\omega$}, Eq.~(\ref{m2h2}),
        as compared to the {\sl quadratic deterioration with~$\omega$\/} of
        the exponentially fitted counterpart.
       }
  \item{The errors associated to the interpolatory Simpson quadrature
        sums~(\ref{qofy}) show quasi-periodic behaviour with respect to
        $\omega$, with a {\sl damping $\omega^{-2}$ power law of the
        amplitudes}, coming from the hypergeometric functions $_0F_1$.
        As a consequence, the accuracy {\sl improves as $\omega$ increases
        over sets of values separated by hypergeometric function induced
        quasi-periods}.

        The best exponentially fitted Simpson quadrature formula $(P\!=\!1)$
        of reference~%
%                     \cite Ixaru [1997]
        \cite{Ixaru97}
        yields comparatively large and roughly uniform error magnitudes over
        sets of $\omega$ values separated by multiples of the $\lambda$-induced
        period $T_{\lambda}$.
       }
\end{enumerate}

Of course, the exponential fitting provides approximating formulas are
significantly better in comparison with classical formulas (wherein the
oscillatory character of the function of interest is ignored). However,
in view of the abovementioned results, the claims made in section~6 of
reference~%
%                     \cite Ixaru [1997]
\cite{Ixaru97}
that the exponentially fitted formulas are "working optimally" on functions
of the form~(\ref{reffunc}) and that further useful extensions can be
developed concerning the numerical differentiation and integration,
are to be regarded with caution.

We end with an interesting observation concerning the interpolatory Simpson
quadrature sums. They can be optimally extended to five-point interpolatory
quadrature sums to yield {\sl quadrature rules\/} along the lines of
{\tt QUADPACK}~%
%                     \cite QUADPACK [1983]
\cite{QUADPACK}.
These can be then conveniently implemented in automatic quadrature codes able
to match the advantages of the polynomial function approximation at
Clenshaw-Curtis or Gauss-Kronrod quadrature knots with the use of the whole
information acquired on the integrand at previous stages of the adaptive
subrange subdivision. This point will be discussed separately in a future paper.
\section{Acknowledgements}
This investigation was partially supported by the Romanian Ministry of
Research and Technology under grant 3036GR/1997.

The graphics was created using the gnuplot package, version 3.7~%
%                       cite{gnuplot}
\cite{gnuplot}.

\end{document}